\documentclass[11pt, a4paper,pagebackref]{amsart}
\usepackage{amsmath}
\usepackage{amssymb}
\usepackage{amsthm}
\usepackage{amscd}
\usepackage{tikz-cd}

\usepackage{color}
\usepackage{hyperref}

\usepackage[inner=2.3cm,outer=2.3cm, bottom=3.3cm]{geometry}
\usepackage{setspace}

\numberwithin{equation}{section}

\usepackage[colorinlistoftodos]{todonotes}

\usepackage[all]{xy}
\usepackage{enumerate}
\usepackage{hyperref}

\usepackage{color}
\hypersetup{colorlinks=true} 

\newtheorem{theorem}{Theorem}[section]

\newtheorem{lemma}[theorem]{Lemma}
\newtheorem{proposition}[theorem]{Proposition}
\newtheorem{corollary}[theorem]{Corollary}

\newcommand{\p}{\mathfrak{p}}

\theoremstyle{definition}
\newtheorem{definition}[theorem]{Definition}
\newtheorem{newclaim}[theorem]{}

\newtheorem{remark}[theorem]{Remark}
\newtheorem{remark and definition}[theorem]{Remark and Definition}
\newtheorem{remark and notation}[theorem]{Remark and Notation}
\newtheorem*{funding}{Funding}
\newtheorem*{agra}{Acknowlegment}

\newtheorem{example}[theorem]{Example}

\newtheorem{question}[theorem]{Question}

\newcommand\Spec{\operatorname{Spec}}

\newcommand\Tor{\operatorname{Tor}}

\newcommand\depth{\operatorname{depth}}

\newcommand\cdim{\operatorname{CI-dim}}

\newcommand\Ker{\operatorname{\Ker}}

\newcommand\pd{\operatorname{pd}}

\DeclareMathOperator{\qpd}{qpd}
\DeclareMathOperator{\tor}{Tor}

\author[Jorge-P\'erez]{Victor H. Jorge-P\'erez}
\author[Martins]{Paulo Martins}
\author[Mendoza-Rubio]{Victor D. Mendoza-Rubio}

\address{Universidade de S{\~a}o Paulo -
ICMC, Caixa Postal 668, 13560-970, S{\~a}o Carlos-SP, Brazil}
\email{vhjperez@icmc.usp.br}

\address{Universidade de S{\~a}o Paulo -
ICMC, Caixa Postal 668, 13560-970, S{\~a}o Carlos-SP, Brazil}
\email{paulomartinsmtm@gmail.com}
\address{Universidade de S{\~a}o Paulo -
ICMC, Caixa Postal 668, 13560-970, S{\~a}o Carlos-SP, Brazil}
\email{vicdamenru@gmail.com}


\begin{document}

\title{Remarks on Auslander's depth formula for quasi-projective dimension}

\maketitle

\begin{abstract}
For nonzero finitely generated $R$-modules $M$ and $N$ over a Noetherian local ring $R$, Auslander's depth formula is the equality
$$
\operatorname{depth} M + \operatorname{depth} N
= \operatorname{depth} R + \operatorname{depth}(\operatorname{Tor}_q^R(M,N)) - q,
$$
where $ q := \sup\{ i \ge 0 \mid \operatorname{Tor}_i^R(M,N) \neq 0 \}$.
Gheibi, Jorgensen, and Takahashi introduced a homological invariant called quasi-projective dimension, which generalizes projective dimension, and proved that Auslander's depth formula holds when $M$ has finite quasi-projective dimension and $q=0$.
In this paper, we prove that the formula still holds when $M$ has finite quasi-projective dimension, $q<\infty$ and $\operatorname{depth}(\operatorname{Tor}_q^R(M,N)) \leq 1$.
We present several applications of this result; in particular, we recover a theorem of Araya and Yoshino, extend our result to the setting of semidualizing modules, and in this framework derive an improved version of the dependency formula for quasi-projective dimension with respect to a semidualizing module recently obtained by Dey, Ferraro, and Gheibi.
\end{abstract}

\section{Introduction}

Throughout this paper, we assume $R$ to be a commutative Noetherian local ring with maximal ideal $\mathfrak{m}$, and that all $R$-modules are finitely generated. Let $M$ and $N$ be $R$-modules, we define: $$q^R(M,N):=\sup \lbrace i \geq 0 : \operatorname{Tor}_i^R(M,N) \neq 0 \rbrace.$$ In one of his renowned papers, Auslander established a classical result that gave rise to the depth formula. This formula plays a crucial role in the study of module depth by relating it to $\Tor$ functors, thus providing an essential understanding of the structure and properties of modules over local rings. Precisely, the theorem is stated as follows:

\begin{theorem}\label{austordep} \cite[Theorem 1.2]{aus} Let $M$ and $N$ be non-zero modules over the local ring $R$ such that $\operatorname{pd}_R M < \infty$. Let $q=q^R(M,N)$. If  $\depth (\tor_q^R(M,N))\le 1$ or $q=0$, then we have the following equality
\begin{align}\label{formula}
 \depth N=\depth(\tor_q^R(M,N)) +\pd_RM-q.   
\end{align}
 \end{theorem}
Note that by Auslander-Buschbaum formula, one can  rewritten \eqref{formula}  as follows:
\begin{align}\label{depthformulag}
\operatorname{depth} N + \operatorname{depth} M = \operatorname{depth} R+ \operatorname{depth}(\operatorname{Tor}_q^R(M,N))-q.
\end{align}

Subsequently, the following question has been studied so far in several articles (see e.g \cite{ArayaYoshino,Celikbas3,Choi, Jorgensen,GJT, HunekeWie,SINHAOthers}). 
\begin{question}\label{q}
    Let $M$ and $N$ be non-zero $R$-modules and let $q=q^R(M,N)<\infty$. Under what conditions, does  the equality
\begin{align}\label{eqdepthfge}
\operatorname{depth} N + \operatorname{depth} M = \operatorname{depth} R+ \operatorname{depth}(\operatorname{Tor}_q^R(M,N)) -q
\end{align}
hold?    
\end{question}
Araya and Yoshino \cite[Theorem 2.5]{ArayaYoshino} provided a positive answer to Question \ref{q} with the conditions that $\operatorname{CI-dim}_R M<\infty$, $q< \infty$ and either $q=0$ or $\operatorname{depth}(\Tor_q^R(M,N))\leq 1$. Bergh and Jorgensen \cite[Theorem 3.1]{Bergh} gave a partial answer to the above question over Cohen-Macaulay local rings with the conditions that $M$ has reducible complexity, $q< \infty$ and either we assume (1) $\operatorname{depth}(\operatorname{Tor}_q^R (M,N))=0$ or (2) $q\geq 1$, $\operatorname{depth}(\operatorname{Tor}_q^R (M,N)) \leq 1$ and $\operatorname{reddeg*} M \geq 2$. Celikbas et al. \cite[Corollary 4.2]{Celikbas3} have generalized the result of Bergh and Jorgensen, proving that Question \ref{q} holds over Cohen-Macaulay local rings, provided that $M$ has finite reducing
complete intersection dimension, $q< \infty$ and $\operatorname{depth} (\operatorname{Tor}_q^R(M,N)) \leq 1$. In the case where  $q=0$, that is $\operatorname{Tor}_{i>0}^R(M,N)=0$, positive results have been established in \cite{Jorgensen,remarks,GJT}, among others. Nevertheless, it worth to  mention that validity of the formula \eqref{eqdepthfge} in the case $q=0$ has been an open problem until that Kimura, Lyle, Otake and Takahashi \cite{kimura2023vanishing} proved that it is not always true.

Recently, Gheibi, Jorgensen, and Takahashi \cite{GJT} introduced a new homological dimension called quasi-projective dimension (qpd) which is a finer invariant than projective dimension, in the sense that for any $R$-module $M$, there is always an inequality $\operatorname{qpd}_R M \leq \operatorname{pd}_R M$ and equality holds when $\operatorname{pd}_R M $ is finite. Among other results, they established a version of the Auslander-Buchsbaum formula for quasi-projective dimension and a version of Auslander's depth formula (\ref{depthformulag}) under the assumptions that one of the modules has finite quasi-projective dimension and that $ q^{R}(M,N)=0 $, i.e., $\operatorname{Tor}_{i>0}^{R}(M,N)=0$.
Unlike the projective dimension, if an $R$-module $M$ has finite quasi-projective dimension, it does not guarantee vanishing of the functor $\Tor_i^R(M,-)$ for all sufficiently large $i$.

Motivated by the extensive research and answers both concerning Question \ref{q}, the main goal of this paper is to present a more complete version of Auslander's depth formula for quasi-projective dimension by including the case where $q^R(M,N) < \infty$ and $\operatorname{depth}(\operatorname{Tor}_{q^R(M,N)}^R(M,N)) \le 1$. Thus, we obtain the following theorem which complement \cite[Theorem 4.11]{GJT}, providing a complete extension of Auslander's original result (Theorem \ref{austordep}) for quasi-projective dimension and a new answer to well studied Question \ref{q}.

\begin{theorem}[See Theorem \ref{auslanderqpd}]\label{teo1.3}
  Let $R$ be a local ring, and let $M$  and $N$ be non-zero $R$-modules such that $ q=q^R(M,N)< \infty$, and $\operatorname{qpd}_R M<\infty$. If $\operatorname{depth}(\operatorname{Tor}_q^R(M,N))  \leq 1$, then we have the following equality: 
  $$\operatorname{depth} N = \operatorname{depth}(\operatorname{Tor}_q^R(M,N)) +\operatorname{qpd}_R M-q.$$
\end{theorem}

The approach used to prove the above theorem is inspired by Auslander's proof of Theorem 1.2 in \cite{aus}. As in the proof given  in \cite{aus} by Auslander, we proceed by induction on $\operatorname{depth} N$, first establishing the base case $\operatorname{depth} N = 0$ (see Theorem \ref{prop:general}). 

For justifying the importance of Theorem \ref{teo1.3}, many applications are provided in this paper. Among them we recover a result of Araya and Yoshino, see Corollary \ref{wef}. We also provide applications of Theorem \ref{teo1.3} to the theory of quasi-projective dimension with respect to a semidualizing module, recently introduced in \cite{ferraro}. Namely, we present a more general version of Auslander's depth formula for quasi-projective dimension with respect to a semidualizing module, complementing the case proved for Tor-independent modules under certain additional conditions in \cite[Theorem 6.10]{ferraro}.

\begin{theorem}[{See Theorem \ref{depthformulaC}}] \label{teo:1.4}
Let $R$ be a local ring, $M$ and $N$ be non-zero $R$-modules such that $q^R(M,N)< \infty$ and let $C$ be a semidualizing $R$-module. If the following conditions are satisfied 
\begin{enumerate}
    \item $\operatorname{\textit{C}-qpd}_R M < \infty$ and $M \in \mathcal{B}_C(R)$, 
    \item $N \in \mathcal{A}_C(R)$, 
    \item $q^R(M,N)=0$ or $\operatorname{depth}(\operatorname{Tor}_{q^R(M,N)}^{R}(M,N)) \leq 1$,
\end{enumerate}
then $$\operatorname{depth} N + \operatorname{depth} M = \operatorname{depth} R+ \operatorname{depth}(\operatorname{Tor}_{q^R(M,N)}^R(M,N)) -q^R(M,N).$$
\end{theorem}

As an application of Theorem \ref{depthformulaC}, we provide a new proof of the dependency formula for $C$-quasi-projective dimension (\cite[Corollary 6.15]{ferraro}). Moreover, our result improves the result of Dey, Ferraro, and Gheibi by showing that $q^R(M,N)$ can be computed from the (finitely many) prime ideals in $\operatorname{Ass}(\operatorname{Tor}_{q^R(M,N)}^R(M,N))$.

\begin{theorem}[{See Theorem \ref{depenency}}] \label{teo1.5} Let $R$ be a local ring, $M$ and $N$ be non-zero $R$-modules such that $q^R(M,N)< \infty$ and let $C$ be a semidualizing $R$-module. If the following conditions are satisfied 
\begin{enumerate}
    \item $\operatorname{\textit{C}-qpd}_R M < \infty$ and $M \in \mathcal{B}_C(R)$, 
    \item $N \in \mathcal{A}_C(R)$,
 \end{enumerate}
 then
 \begin{align*}
    q^R(M,N) & = \sup \lbrace \operatorname{depth} R_{\mathfrak{p}} - \operatorname{depth} M_{\mathfrak{p}} - \operatorname{depth} N_{\mathfrak{p}} \mid \mathfrak{p} \in \Spec R \rbrace \textnormal{ (\cite[Corollary 6.15]{ferraro})} \\
    & =  \sup \lbrace \operatorname{depth} R_{\mathfrak{p}} - \operatorname{depth} M_{\mathfrak{p}} - \operatorname{depth} N_{\mathfrak{p}} \mid \mathfrak{p} \in \operatorname{Ass} (\operatorname{Tor}_{q^R(M,N)}^R(M,N)) \rbrace.
    \end{align*}
\end{theorem}

    Now, we describe the structure of this paper. In Section 2, we establish conventions and provide some facts about the quasi-projective dimension that are used in this paper. In Section 3, we prove Theorem \ref{teo1.3} after first proving some auxiliary results, and we explore some of its consequences. In the last section, we apply Theorem \ref{teo1.3}  to theory of quasi-projective dimension with respect to a semidualizing module; specifically proving Theorems \ref{teo:1.4} and \ref{teo1.5}.

\section{Conventions and Background}
In this section, we provide some notation and known facts about quasi-projective dimension, which will play a crucial role in many of our proofs.

\begin{newclaim}
Let $R$ be a local ring and let $M$ be an $R$-module. Consider a minimal free resolution
\begin{align*}
\cdots \rightarrow F_i \xrightarrow{\varphi_i} F_{i-1} \rightarrow \cdots \rightarrow F_1 \xrightarrow{\varphi_1} F_0 \xrightarrow{\varphi_0} M \rightarrow 0
\end{align*}
of $M$. For $i \geq 1$, the $i$-\textit{syzygy} of $M$, denoted by $\Omega^i(M)$, is defined as the kernel of the map $\varphi_{i-1}$. When $i=0$, we set $\Omega^0(M)=M$. For $i \geq 0$, the modules $\Omega^i(M)$ are defined uniquely up to isomorphism.
\end{newclaim}

\begin{newclaim}

For a complex $$X: \cdots \stackrel{\partial_{i+2}}{\longrightarrow} X_{i+1} \stackrel{\partial_{i+1}}{\longrightarrow  } X_i \stackrel{\partial_{i}}{\longrightarrow
      } X_{i-1}   \longrightarrow  \cdots $$ of $R$-modules, we set for each integer $i$, $\operatorname{Z}_i(X)=\ker \partial_i$ and $\operatorname{B}_i(X)= \operatorname{Im} \partial_{i+1}$ and $\operatorname{H}_i(X)=\operatorname{Z}_i(X)/\operatorname{B}_i(X)$.  Moreover, we set: \begin{align*}
&\left\{ 
    \begin{aligned}
        \sup X &= \sup \{ i \in \mathbb{Z} : X_i \neq 0 \},\\
        \inf X &= \inf \{ i \in \mathbb{Z} : X_i \neq 0 \},
    \end{aligned}
\right. \quad 
\left\{
    \begin{aligned}
        \text{hsup } X &= \sup \{ i \in \mathbb{Z} : \operatorname{H}_i(X) \neq 0 \},\\
        \text{hinf } X &= \inf \{ i \in \mathbb{Z} : \operatorname{H}_i(X) \neq 0 \}.
    \end{aligned}
\right.
\end{align*}
The \textit{length} of $X$ is defined to be  $\operatorname{length} X= \sup X - \inf X$. We say that $X$ is \textit{bounded}, if $\operatorname{length} X < \infty$. We say that $X$ is \textit{bounded below} if $\inf X > -\infty$.
\end{newclaim}

\subsection{Quasi-projective dimension.} The definition of quasi-projective dimension was introduced by Gheibi, Jorgensen and Takahashi \cite{GJT} as a generalization of the classical notion of projective dimension.
\begin{definition}(\cite[Definition 3.1]{GJT})\label{defqpd}
    Let $M$ be an $R$-module.
    \begin{enumerate}
        \item A \textit{quasi-projective resolution} of $M$ is bounded below complex $$P: \cdots \longrightarrow P_{i+1} \longrightarrow P_i \longrightarrow P_{i-1} \longrightarrow \cdots $$ of projective $R$-modules such that for all $i\geq \inf P$, there exist non-negative integers $a_i$, not all zero, such that $\operatorname{H}_i(P) \cong M^{\oplus a_i}$. 
        
    \item      The \textit{quasi-projective dimension} of $M$, denoted by $\operatorname{qpd}_R M$ is defined as follows: 
        $$\mathrm{qpd}_R M = \inf \{ \sup P - \text{hsup } P \mid P \text{ is a bounded quasi-projective resolution of } M \}$$  if  $M\not=0$, and $\operatorname{qpd}_R M=-\infty$ if $M=0$.
    \end{enumerate}
\end{definition}
One has $\text{qpd}_R M = \infty$ if and only if $M$ does not admit a bounded quasi-projective resolution.
In \cite{GJT}, the authors established several results, including versions of the Auslander-Buchsbaum formula, Auslander's depth formula, and the Auslander-Reiten conjecture for modules of finite quasi-projective dimension. Some of these results are stated below.

A complex $(X,\partial)$ of free $R$-modules of finite rank is called \textit{minimal} if $\partial_i \otimes_R k=0$ for all $i$, where $k$ is the residue field of $R$. 

\begin{proposition} \label{prop41}\cite[Proposition 4.1]{GJT}
Let $R$ be a local ring, and let $M$ be a non-zero $R$-module with $\mathrm{qpd}_R M < \infty$. Then there exists a finite minimal quasi-projective resolution $F$ of $M$ such that $\mathrm{qpd}_R M = \sup F - \mathrm{hsup} F$.
\end{proposition}

\begin{theorem}\cite[Theorem 4.4]{GJT}\label{rem:ABF}
Let $R$ be a local ring, and let $M$ be an $R$-module of finite quasi-projective dimension. Then
\[
\operatorname{qpd}_R M = \operatorname{depth} R - \operatorname{depth} M.
\]
In particular, if $M\not=0$, then $\operatorname{depth} M\leq \operatorname{depth} R$ and $\operatorname{qpd}_R M\leq \operatorname{depth} R.$
\end{theorem}

\begin{theorem}\cite[Theorem 4.11]{GJT}\label{depthformula}
Let $R$ be a local ring, and let $M$ and $N$ be $R$-modules. Suppose that $M$ has finite quasi-projective dimension and $\operatorname{Tor}_{i}^R(M,N)=0$ for all $i>0$. Then 
\[
\operatorname{depth} M + \operatorname{depth} N = \operatorname{depth} R + \operatorname{depth} (M \otimes_R N).
\]
\end{theorem}

\section{Auslander's depth formula for modules of finite quasi-projective dimension}  

In this section, we prove Theorem \ref{teo1.3}, which is Theorem \ref{auslanderqpd} in this section, and explore some of its consequences. Before proving  Theorem \ref{teo1.3}, we prove the following theorem,  which show Theorem \ref{teo1.3} holds in the case where $\depth N=0$ and will be helpful in its proof. The next theorem is a version of Auslander's result \cite[Proposition 1.1]{aus} for quasi-projective dimension.

\begin{theorem}\label{prop:general}
Let $R$ be a local ring, and let $M$ and $N$ be non-zero $R$-modules such that $q^R(M,N) < \infty$ and $\operatorname{depth} N=0$. Suppose that $r=\operatorname{qpd}_R M < \infty$. Then $\operatorname{Tor}_{r}^R(M,N)\not=0$ and has depth zero. 
\end{theorem}

For proving Theorem \ref{prop:general}, we need establishing some auxiliary results.

\begin{lemma}\label{sel}
    Let $R$ be a local ring and let $M$ and $N$ be $R$-modules such that $q^R(M,N)<\infty$. Then $q^R(M,N)\leq \operatorname{qpd}_R M$.
\end{lemma}
\begin{proof}
We may assume that $\operatorname{qpd}_R M< \infty$. Since $\operatorname{Tor}_{i\gg0}^R(M,N)=0$, then by \cite[Corollary 6.4]{GJT}, we have that $\operatorname{Tor}_i^R(M,N)=0$ for all $i \geq \operatorname{qpd}_R M+1$ and the desired inequality follows.
\end{proof}

\begin{proposition}\label{analemnaus}
    Let $R$ be a local ring, and let $M$ and $N$ be non-zero $R$-modules such that $q^R(M,N)<\infty$. Assume $\depth N=0$. If 
$\operatorname{qpd}_R M<\infty$, then $q^R(M,N)=\operatorname{qpd}_R M$.  
\end{proposition}
\begin{proof}
Since $M$ has finite quasi-projective dimension, by Proposition \ref{prop41}, there exists a minimal quasi-projective resolution $(F, \partial)$ of $M$ with $s=\sup F,$ $h=\operatorname{hsup} F$ and $r:=\operatorname{qpd}_R M=s-h$. Let $C=\operatorname{Coker}(\partial_{h+1})$ and $q=q^R(M,N)$. Then we obtain an exact sequence 
$$0 \rightarrow F_s \xrightarrow{\partial_s} \cdots \xrightarrow{\partial_{h+1}} F_h \rightarrow C \rightarrow 0,$$
which is a minimal free resolution of $C$, so that $\operatorname{pd}_R C=s-h=r$.
For all $i$, set $Z_i=Z_i(F)$ and $B_i=B_i(F)$. There are exact sequences
\begin{align*}
    \begin{cases}
0 \rightarrow B_i \rightarrow Z_i \rightarrow M^{\oplus a_i} \rightarrow 0 \\  0 \rightarrow Z_i \rightarrow F_i \rightarrow B_{i-1} \rightarrow 0
\end{cases}  
\quad (i \in \mathbb{Z})
\end{align*}
with $a_i \geq 0$. Since $\operatorname{Tor}_i^R(M,N)=0$ for all $i>q$, we see that 
\[
\operatorname{Tor}_j^R(B_i, N) \cong \operatorname{Tor}_j^R(Z_i, N) \cong \operatorname{Tor}_{j+1}^R(B_{i-1}, N) \cong \cdots \cong \operatorname{Tor}_{j+i+1-\inf F}^R(B_{\inf F - 1}, N) = 0
\]
for all \( i \geq \inf F \) and \( j > q \), where the last equality holds because \( B_{\inf F-1} = 0 \). Thus, $\operatorname{Tor}_j^R(B_i, N)=0$ for all $i$ and $j>q$.

Note that we also have an exact sequence
\begin{equation*} 
    0 \rightarrow \operatorname{H}_h(F) \rightarrow C \rightarrow \operatorname{B}_{h-1} \rightarrow 0.
\end{equation*} Since $\operatorname{Tor}_j^R(B_{h-1}, N)=0=\Tor_j^R(M,N)$ for all $j>q$, then $\Tor_j^R(C,N)=0$ for all $j>q$. Thus, if $r>q$, then $\Tor_r^R(C,N)=0$, which  contradicts to \cite[Proposition 1.1]{aus} as $\operatorname{pd}_R C=r$. Therefore $r\leq q$ and it follows from Lemma \ref{sel} that $r=q$.
\end{proof}

The following is a consequence of Proposition \ref{analemnaus},  and  although it is not used in the proof of our main theorem, it will be used in one of its applications.
\begin{corollary} \label{ineqgs} Let $R$ be a local ring, and let $M$ and $N$ be non-zero $R$-modules such that $q^R(M,N)<\infty$.  If 
$\operatorname{qpd}_R M<\infty$, then $q^R(M,N) \geq \operatorname{qpd}_R M - \operatorname{depth} N$.
\end{corollary}
\begin{proof}
    If $\depth N=0$, the result follows directly from the previous proposition. Now, suppose that $\operatorname{depth} N >0$. Let $x \in \mathfrak{m}$ be an $N$-regular element. Consider the following exact sequence $0 \rightarrow N \xrightarrow{.x}  N \rightarrow N/xN \rightarrow 0 $. From its long exact sequence $\operatorname{Tor}^R(M,-)$, it is easy to see that $q^R(M,N)+ 1 \geq q^R(M,N/xN)$. Then, by induction we have
    \begin{align*}
    q^R(M,N)+1 \geq q^R(M,N/xN) \geq \operatorname{qpd}_R M- \operatorname{depth} N/xN.
\end{align*}
Thus $q^R(M,N) \geq \operatorname{qpd}_R M - \operatorname{depth} N$.

\end{proof}

\begin{lemma} \label{lemavhcorr}
Let \( R \) be a local ring, and let \( M \) and \( N \) be non-zero \( R \)-modules such that \(\operatorname{depth} N=0\) and \(\qpd_R M < \infty\). If $q^R(M,N)=1$, then $\operatorname{depth} (\operatorname{Tor}_1^R(M,N))=0$.
\end{lemma}
\begin{proof}
    By Proposition  \ref{analemnaus},  $\operatorname{qpd}_R M=1$. Then there exists a minimal quasi-projective resolution $(F,\partial)$ such that $1=\qpd_R M=\operatorname{sup} F-\operatorname{hsup} F$. Set $s=\sup F$ and $h=\operatorname{hsup}$. Let  $C=\operatorname{coker}(\partial_{h+1})$. Then we have an exact sequence 
    $$0 \rightarrow F_s \xrightarrow{\partial_s} F_h \rightarrow C \rightarrow 0.$$
   Again, since $F$ is minimal, we see that $\operatorname{pd}_R C=1$. By \cite[Lemma 2.5]{GJT}, there exists a convergent spectral sequence

    $$\operatorname{E}_{i,j}^2=\operatorname{Tor}_i^R(\operatorname{H}_j(F), N) \Rightarrow \operatorname{H}_{i+j}(F \otimes_R N).$$

    Since each $\operatorname{H}_j(F)$ is a finite direct sum of copies of $M$, by assumption, we have that $\operatorname{E}^2_{i,j}=0$ for $i\not=0,1$. Therefore, the spectral sequence induces, for each $n$, an exact sequence
    $$0 \rightarrow \operatorname{H}_n(F) \otimes_R N
  \rightarrow \operatorname{H}_n(F \otimes_R N) \rightarrow  \Tor_1^R(\operatorname{H}_{n-1}(F), N) \rightarrow 0 .$$
  In particular, if $n=s=h+1$, then the first term is zero and we have isomorphisms
  $$\operatorname{H}_s(F \otimes_R N) \cong \operatorname{Tor}_1^R(\operatorname{H}_{h}(F),N )\cong \Tor_1^R(M^{\oplus a_h}, N)\cong \operatorname{Tor}_1^R(M,N)^{\oplus a_h},$$
for a positive integer $a_h$. Also,  note that $$\operatorname{H}_s(F \otimes_R N)=\operatorname{ker}(\partial_s \otimes_R \operatorname{Id}_N)\cong \operatorname{Tor}_1^R(C,N).$$ Thus
  $\operatorname{Tor}_1^R(C,N)\cong \operatorname{Tor}_1^R(M,N)^{\oplus a_h}$, for $a_h>0$. By \cite[Proposition 1.1]{aus}, the module $\operatorname{Tor}_1^R(C,N)$ has depth zero, and hence $\operatorname{Tor}_1^R(M,N)$ as well. 
\end{proof}

\begin{proof}[Proof of Theorem \ref{prop:general}]
Let $q=q^R(M,N)$. By Proposition \ref{analemnaus}, we have $q=\operatorname{qpd}_R M$. We prove the theorem by considering the following two cases:

\begin{enumerate}
\item If $q=0,$ then  $\operatorname{Tor}_i^R(M,N)=0$ for all $i>0$. Therefore, we see from the Auslander-Buschbaum formula for quasi-projective dimension that $\operatorname{depth} R=\operatorname{depth} M$. Thus, it follows from Theorem \ref{depthformula} that $\operatorname{Tor}_q^R(M,N)\cong M \otimes_R N$ has depth zero.
 
  \item Suppose $q>0$. Note that  $\Tor_1^R(\Omega^{q-1} M, N)\cong \Tor_q^R(M,N)\not=0$ and that $\Tor_i^R(\Omega^{q-1} M, N) \cong \Tor_{i+q-1}^R(M,N)=0$ for $i\geq 2$. Thus $\Omega^{q-1} M \not=0$ and $q^R(\Omega^{q-1} M,N)=1$. By \cite[Corollary 4.5]{GJT}, we see that $\operatorname{qpd}_R(\Omega^{q-1} M)=1<\infty$. Then, by Lemma \ref{lemavhcorr}, we have that $ \Tor_q^R(M,N) \cong \Tor_1^R(\Omega^{q-1} M, N)$ has depth zero. 
\end{enumerate}    
\end{proof}

\begin{remark}\label{victorremark}
Let $R$ be a local ring and let $x \in \mathfrak{m}$ a non-zero element of its maximal ideal. Let $M$ be an $R$-module of depth zero. Consider the exact sequence 
\begin{align*}
0 \rightarrow K \rightarrow M \xrightarrow{\cdot x} M
\end{align*}
where $K$ is the kernel of the multiplication by $x$ on $M$. Then $K \neq 0$ has depth zero. Indeed, by contradiction, assume that $\operatorname{depth} K >0$. As $\operatorname{depth} M =0$, then $\mathfrak{m}=\operatorname{Ann}(m)$, for a non-zero element $m \in M$. As $x \in \mathfrak{m}$, then $xm=0$ and  hence $m \in K$. Also, as $\operatorname{depth}K>0$, we can find an element $y \in \mathfrak{m}$ that is $K$-regular. Note that $ym=0$ and it is a contradiction.
\end{remark}
\begin{theorem}\label{auslanderqpd}
Let $R$ be a local ring, and let $M$  and $N$ be non-zero $R$-modules such that $q=q^R(M,N)<\infty$ and $\operatorname{qpd}_R M<\infty$. If $\operatorname{depth} (\operatorname{Tor}_q^R(M,N))  \leq 1$, then we have the following equality: 
  $$\operatorname{depth} N = \operatorname{depth}(\operatorname{Tor}_q^R(M,N)) +\operatorname{qpd}_R M-q.$$
\end{theorem}
\begin{proof}
 Assume that $\operatorname{depth} N=0$. From Proposition \ref{analemnaus} and Theorem \ref{prop:general}, we see that $q=\operatorname{qpd}_R M$ and that $\operatorname{Tor}_q^R(M,N)$ has depth zero. Therefore, the equality holds. Suppose that $\operatorname{depth} N = k > 0$. Let $x \in \mathfrak{m}$ be an $N$-regular element. The exact sequence 
\begin{align*}
0 \rightarrow N \xrightarrow{\cdot x} N \rightarrow N/xN \rightarrow 0
\end{align*}
induces the exact sequence 
\begin{align}\label{seq:exata}
\cdots \rightarrow  \operatorname{Tor}_{q+1}^R(M,N) \rightarrow \operatorname{Tor}_{q+1}^R(M,N/xN) \rightarrow \operatorname{Tor}_q^R(M,N) 
\end{align}
\begin{align*}
 \quad \quad    \xrightarrow{\cdot x} \operatorname{Tor}_q^R(M,N) \rightarrow 
  \operatorname{Tor}_q^R(M,N/xN) \rightarrow \cdots.    
\end{align*}
Therefore, we have that $\operatorname{Tor}_i^R(M,N/xN)=0$ for all $i > q+1$, since $\operatorname{Tor}_i^R(M,N)=0$, for $i>q$. We will consider separately the cases $\operatorname{depth}(\operatorname{Tor}_q^R(M,N))=0$ and $\operatorname{depth}(\operatorname{Tor}_q^R(M,N))=1$. 

In the first case, we proceed by induction on $\operatorname{depth} N$. If $\operatorname{depth}( \operatorname{Tor}_q^R(M,N))=0$, then $\operatorname{Tor}_{q+1}^R(M,N/xN) \neq 0$ and has depth zero (see Remark \ref{victorremark}). Since $\operatorname{depth}(N/xN)=\operatorname{depth} N -1$, then applying induction we obtain
$$\operatorname{depth} N -1 = \operatorname{depth} (\operatorname{Tor}_{q+1}^R(M,N/xN)) + \operatorname{qpd}_R M -(q+1),$$
which implies the desired equality.

Now, if $\operatorname{depth}(\operatorname{Tor}_q^R(M,N))=1$, we may assume that $x$ is also regular on  $\operatorname{Tor}_q^R(M,N)$. From the long exact sequence \eqref{seq:exata}, we obtain that $\operatorname{Tor}_i^R(M,N/xN)=0$ for $i > q$. Considering again the exact sequence \eqref{seq:exata}, one can conclude using Nakayama's Lemma that $\operatorname{Tor}_q^R(M,N/xN) \neq 0$ and obtain the exact sequence
\begin{align}\label{seque:exata2}
0 \rightarrow \operatorname{Tor}_q^R(M,N) \xrightarrow{\cdot x} \operatorname{Tor}_q^R(M,N) \rightarrow \operatorname{Tor}_q^R(M,N/xN).
\end{align}
It follows from the exact sequence \eqref{seque:exata2} that $\operatorname{Tor}_q^R(M,N/xN)$ contains a submodule of depth zero and thus $\operatorname{depth}(\operatorname{Tor}_q^R(M,N/xN))=0$. By the previous case, we have that
\begin{align*}
\operatorname{depth} N -1 = \operatorname{depth}(\operatorname{Tor}_q^R(M,N/xN)) + \operatorname{qpd}_R M -q
\end{align*}
and
    $$\operatorname{depth} N = \operatorname{depth}(\operatorname{Tor}_q^R(M,N)) +\operatorname{qpd}_R M-q.$$
\end{proof}

The following example demonstrates that the condition $q^R(M,N)< \infty$ in Theorem \ref{auslanderqpd} cannot be omitted.

\begin{example}
Let $(R,\mathfrak{m},k)$ be a non-regular local ring. The residue field $k$ has finite quasi-projective dimension \cite[Proposition 3.6(1)]{GJT}. On the other hand, since $R$ is non-regular, we must have that $\operatorname{pd}_R k= \infty$ and $q^R(k,k)= \infty$. 
\end{example}

Now, we provide some consequences of Theorem \ref{auslanderqpd}. The next corollary extends \cite[Theorem 4.11]{GJT} and provides a new answer to Question \ref{q}.
\begin{corollary}\label{wef}
      Let $R$ be a local ring, and let $M$  and $N$ be non-zero $R$-modules such that $ q=q^R(M,N)< \infty$ and $\operatorname{qpd}_R M<\infty$. If  $\operatorname{depth} (\operatorname{Tor}_q^R(M,N))  \leq 1$ or $q=0$, then we have the following equality:  $$\operatorname{depth} N + \operatorname{depth} M = \operatorname{depth} R+ \operatorname{depth}(\operatorname{Tor}_q^R(M,N)) -q.$$
\end{corollary}
\begin{proof} 
The case $q=0$ is Theorem \ref{depthformula} and the case $\operatorname{depth} (\operatorname{Tor}_q^R(M,N))  \leq 1$ follows from Theorems \ref{rem:ABF} and \ref{auslanderqpd}.  
\end{proof}

We say that a ring $S$ is a \textit{deformation} of a ring $R$ of codimension $r$ if there exists an $S$-regular sequence $\textbf{x}=x_1,\dots,x_r$ such that $S/(\textbf{x}) \cong R$. And we say that $R \rightarrow R' \leftarrow S$ is a \textit{quasi-deformation} of codimension $r$ if $R \rightarrow R'$ is a local flat homomorphism and $S$ is a deformation of $R'$ of codimension $r$. The complete intersection dimension of an $R$-module $M$ is defined as follows: 
\begin{align*}
\cdim_R M = \inf \lbrace \pd_S (M \otimes_R R')- \operatorname{pd}_S R' \mid R \rightarrow R' \leftarrow S \text{ is a quasi deformation}\rbrace.
\end{align*}

Araya and Yoshino \cite[Theorem 2.5]{ArayaYoshino} proved the same above result with complete intersection dimension instead of quasi-projective dimension. As an application of Corollary \ref{wef}, we recover this result.
\begin{corollary}\label{cor:araya}\cite[Theorem 2.5]{ArayaYoshino} 
    Let $R$ be a local ring, and let $M$  and $N$ be non-zero $R$-modules such that $ q=q^R(M,N)<\infty$ and  $\cdim_R M<\infty$. If  $\operatorname{depth} (\operatorname{Tor}_q^R(M,N))  \leq 1$ or $q=0$, then we have the following equality:  $$\operatorname{depth} N + \operatorname{depth} M = \operatorname{depth} R+ \operatorname{depth}(\operatorname{Tor}_q^R(M,N)) -q.$$
\end{corollary}
\begin{proof}
Since $\cdim_R M < \infty$, there exists a quasi-deformation $R \rightarrow R' \leftarrow S$ such that $\operatorname{pd}_S (M \otimes_R R')< \infty$. As at the beginning of the proof of \cite[Theorem 3]{Choi}, without loss of generality, we may assume that $R=R'$, i.e, $S$ is a deformation of $R$. Since $\operatorname{pd}_S M < \infty$, then $\operatorname{qpd}_RM < \infty$ by \cite[Proposition 3.7]{GJT}. Thus, the desired formula follows directly from Corollary \ref{wef}.
\end{proof}

\begin{corollary}\label{corol:intersection}
Let $R$ be a local ring, and let $M$  and $N$ be non-zero $R$-modules such that $ q=q^R(M,N)<\infty$ and $\operatorname{qpd}_R M<\infty$. Assume that $q \geq 1$. If $\operatorname{depth} (\operatorname{Tor}_q^R(M,N))  \leq 1$, then $\operatorname{depth} N \leq \operatorname{qpd}_R M$. In particular, if $N$ is Cohen-Macaulay, then $\dim N \leq \operatorname{qpd}_R M$.
\begin{proof}
Follows by Theorem \ref{auslanderqpd}.
\end{proof}
\end{corollary}
Corollary \ref{corol:intersection} leads us to consider the following question:
\begin{question}\label{question}
Under what conditions can we assert that $\dim N \leq \operatorname{qpd}_R M$, as stated in the corollary above, or that $\dim N \leq \operatorname{qpd}_R M + \dim (M \otimes_R N)$?
\end{question}
\begin{example}
It is easy to see that the inequalities of Question \ref{question} are not satisfied in general. For example:
\begin{enumerate}
    \item The residue field $k$ of a local ring $R$ always has finite quasi-projective dimension (see \cite[Proposition 3.6]{GJT}). The inequality $\dim R \leq \operatorname{qpd}_R k=\operatorname{depth} R$ is not satisfied unless $R$ is Cohen-Macaulay.
    \item Let $k$ be a field and let $R=k[[X,Y]]/(XY)$, $M=R/(X)$ and $N=R/(Y)$. By \cite[Corollary 3.8]{GJT} $\operatorname{qpd}_R M < \infty$. It follows by the Auslander-Buchsbaum formula for quasi-projective dimension that $\operatorname{qpd}_R M =0$. Also, it is easy to see that $\dim N =1$ and $\operatorname{dim} M \otimes_R N =0$. Then, we have that $\dim N > \operatorname{qpd}_R M + \dim M \otimes_R N$. 
\end{enumerate}
\end{example}
\begin{remark}
A partial answer to Question \ref{question} was recently provided by the authors in \cite[Corollary 3.5]{VPV}.    
\end{remark}

Let $M$ and $N$ be $R$-modules. If $\operatorname{Tor}_{i>0}^R(M,N)=0$, then we say that $M$ and $N$ are $\operatorname{Tor}$\textit{-independent}. Moreover, an $R$-module $M$ satisfies \textit{Serre's} $(S_n)$\textit{-condition} for a positive integer $n$ provided $\depth M_{\p} \geq \min \lbrace n, \dim R_{\p} \rbrace$, for all $\p \in \Spec (R)$.  
\begin{corollary}
Let $R$ be a Cohen-Macaulay local ring, and let $M$ and $N$ be $R$-modules. Assume that $\operatorname{qpd}_R M < \infty$ and that $M$ and $N$ are $\operatorname{Tor}$-independent. Then if $M \otimes_R N$ satisfies Serre's $(S_n)$-condition for some $n$, so do the modules $M$ and $N$. 
\end{corollary}
\begin{proof}
For any $\mathfrak{p} \in \operatorname{Supp}(M \otimes_R N)$, we have that $\operatorname{qpd}_{R_{\mathfrak{p}}}M_{\mathfrak{p}} \leq \operatorname{qpd}_R M < \infty$ by \cite[Proposition 3.5]{GJT}. Also, it is easy to see that $M_{\mathfrak{p}}$ and $N_{\mathfrak{p}}$ are $\operatorname{Tor}$-independent over $R_{\mathfrak{p}}$. Hence, by Corollary \ref{wef}, we have 
\begin{align*}
\operatorname{depth} M_{\mathfrak{p}}+\operatorname{depth} N_{\mathfrak{p}} =\operatorname{depth} R_{\mathfrak{p}}+ \operatorname{depth} (M \otimes_R N)_{\mathfrak{p}} \geq \operatorname{depth} R_{\mathfrak{p}}+ \min \lbrace n, \dim R_{\mathfrak{p}} \rbrace.
\end{align*}
Since $R$ is Cohen-Macaulay, we can see that $R_{\mathfrak{p}}$ is Cohen-Macaulay and consequently we have $\operatorname{depth} M_{\mathfrak{p}},\operatorname{depth} N_{\mathfrak{p}} \leq \operatorname{depth} R_{\mathfrak{p}} $. Then, we have $\operatorname{depth} M_{\mathfrak{p}},$ $\operatorname{depth} N_{\mathfrak{p}} \geq \min \lbrace n, \dim R_{\mathfrak{p}} \rbrace$. 
\end{proof}
Recall that an $R$-module $M$ is \textit{torsion} if  $M$ is equals its torsion submodule, which is equivalent to the fact that  $M_{\mathfrak{p}}=0$ for all $\mathfrak{p} \in \operatorname{Ass}(R)$. The next result is comparable to \cite[A2]{Celikbas2}.
\begin{lemma}\label{torsion}
Let $R$ be a local ring and let $M$ and $N$ be non-zero $R$-modules with $\operatorname{qpd}_R M< \infty$. If $\operatorname{Tor}_i^R(M,N)$ is torsion for all $i \gg 0$, then $\operatorname{Tor}_i^R(M,N) $ is torsion for all $i >0$. 
\end{lemma}
\begin{proof}
Given $\mathfrak{p} \in \operatorname{Ass}(R)$, we shall show that $\Tor_i^{R_\mathfrak{p}}(M_\mathfrak{p},N_\mathfrak{p})=0$ for all $i>0$. We may assume that $M_\mathfrak{p}$ and $N_\mathfrak{p}$ are non-zero.  Considering that $\operatorname{Tor}_i^R(M,N)$ is torsion for all $i \gg 0$, then $q^{R_{\mathfrak{p}}} (M_{\mathfrak{p}},N_{\mathfrak{p}}) < \infty$. By \cite[Proposition 3.5(1)]{GJT}, we also have $\operatorname{qpd}_{R_{\mathfrak{p}}}M_{\mathfrak{p}} < \infty$. Using Lemma \ref{sel} and Theorem \ref{rem:ABF}, we have that $q^{R_{\mathfrak{p}}} (M_{\mathfrak{p}},N_{\mathfrak{p}}) \leq \operatorname{depth} R_{\mathfrak{p}}=0$ and hence $q^{R_{\mathfrak{p}}} (M_{\mathfrak{p}},N_{\mathfrak{p}})=0$. 
\end{proof}

\begin{proposition}
Let $R$ be a one-dimensional Cohen-Macaulay local ring and let $M$ and $N$ be non-zero $R$-modules. If $q^R(M,N)<\infty$ and $\operatorname{qpd}_RM < \infty$, then the following are equivalent:
\begin{enumerate}
    \item $q^R(M,N)=0$.
    \item $\operatorname{depth} M + \operatorname{depth} N = \operatorname{depth} R + \operatorname{depth} M \otimes_R N$.
    \item  $M$ or $N$ is torsion-free.
\end{enumerate}
\end{proposition}
\begin{proof}
The assertion (1) $\Rightarrow$ (2) follows by the depth formula for quasi-projective dimension (see Corollary \ref{wef} or Theorem \ref{depthformula}). Now assume that (2) holds true. In this case, we have $2 \geq \operatorname{depth} M + \operatorname{depth} N \geq \operatorname{depth} R =1$ and thus $\operatorname{depth} M =1$ or $\operatorname{depth} N =1$ and the statement of (3) is verified. 

To end this proof, assume (3) holds true. By Lemma \ref{sel} and Theorem \ref{rem:ABF}, we have that $q^R(M,N) \leq \operatorname{qpd}_R M \leq 1$. By contradiction, assume that $q^R(M,N) \neq 0$, that is, $q^R(M,N)=1$. As $\operatorname{depth} (\operatorname{Tor}_1^R(M,N)) \leq 1$, our generalization of depth formula for quasi-projective dimension (Theorem \ref{auslanderqpd}) provides the equality
$$
\operatorname{depth} M + \operatorname{depth} N = \operatorname{depth} (\operatorname{Tor}_1^R(M,N)).
$$
By Lemma \ref{torsion}, we have that $\operatorname{Tor}_1^R(M,N)$ is torsion and then has finite length and depth zero. Thus $\operatorname{depth} M + \operatorname{depth} N =0$ and it is a contradiction, since $M$ or $N$ is torsion-free. That is, the implication (3) $\Rightarrow (1)$ is proved.
\end{proof}
The next theorem is a version of \cite[Theorem 2.10]{ArayaYoshino} with quasi-projective dimension instead of complete intersection dimension. Our Theorem \ref{auslanderqpd} plays a crucial role on its proof. Recall that a local ring $(R,\mathfrak{m})$ is an \textit{isolated singularity} provided $R_{\mathfrak{p}}$ is a regular local ring for every prime ideal $\mathfrak{p} $ in the punctured spectrum $\operatorname{Spec} R \backslash \lbrace \mathfrak{m} \rbrace$.
\begin{theorem}\label{teo:torindependent}
Let $(R,\mathfrak{m})$ be an isolated singularity, and let $M$ and $N$ be $R$-modules such that $\operatorname{qpd}_R M < \infty$ and $ q^R(M,N)< \infty$. If $\operatorname{depth} M + \operatorname{depth} N \geq \operatorname{depth} R$ and if $M \otimes_R N$ is torsion free on the punctured spectrum, then $M$ and $N$ are $\operatorname{Tor}$-independent.
\end{theorem}
\begin{proof}
Since $R$ is an isolated singularity, then $R_{\mathfrak{p}}$ is a regular local ring, for all prime ideal $\mathfrak{p} \in \operatorname{Spec} R \backslash \lbrace \mathfrak{m} \rbrace$. Also, as $M_{\mathfrak{p}} \otimes_{R_{\mathfrak{p}}} N_{\mathfrak{p}}$ is torsion free, then $M_{\mathfrak{p}}$ and $N_{\mathfrak{p}}$ are $\operatorname{Tor}$-independent over $R_{\mathfrak{p}}$, by \cite[Lemma 2.11]{ArayaYoshino}, i.e, $\operatorname{Tor}_i^R(M,N)_{\mathfrak{p}}\cong \operatorname{Tor}_i^{R_{\mathfrak{p}}}(M_{\mathfrak{p}},N_{\mathfrak{p}})=0$, for all $i>0$. Suppose that $M$ and $N$ are not $\operatorname{Tor}$-independent and denote $q=q^R(M,N)>0$. Since $\Tor_q^R(M,N)$ has finite length, it is clear that $\operatorname{depth} (\Tor_q^R(M,N))=0$. Therefore, we can apply Theorem \ref{auslanderqpd} and thus obtain
\begin{align*}
\operatorname{depth} M + \operatorname{depth} N = \operatorname{depth} R - q < \operatorname{depth} R, 
\end{align*}
which is a contradiction.
\end{proof}
The following corollary is a variation of \cite[Theorem 2.4]{HunekeWie} and \cite[Corollary 2.13]{ArayaYoshino} with quasi-projective dimension instead of projective dimension and complete intersection dimension, respectively. This result is a direct consequence of Corollary \ref{wef} and Theorem \ref{teo:torindependent} and its proof is similar to that in \cite[Corollary 2.13]{ArayaYoshino}.
\begin{corollary}
Let $(R, \mathfrak{m})$ be an isolated singularity, and let $M$ and $N$ be $R$-modules such that $\operatorname{qpd}_R M < \infty$ and $q^R(M,N)<\infty$. Suppose that $\operatorname{depth} M + \operatorname{depth} N \geq \operatorname{depth} R + n$ and that $M \otimes_R N$ satisfies $(S_{n+1})$-condition on the punctured spectrum for some integer $n$ with $0 \leq n < \operatorname{dim} R$. Then the following conditions are equivalent: 
\begin{enumerate}
\item[(1)] $H^n_{\mathfrak{m}} (M \otimes_R N)=0$,
\item[(2)] $M \otimes_R N$ satisfies $(S_{n+1})$-condition,
\item[(3)] $\operatorname{depth}(M \otimes_R N) \geq n+1$. 
\end{enumerate}
\end{corollary}
\begin{proof}
The assertion (2) $\Rightarrow$ (3) follows directly. The implication (3) $\Rightarrow$ (1) is clear, since $H^i_{\mathfrak{m}}(M \otimes_R N)=0$ for $i< \operatorname{depth}(M \otimes_R N)$. We need only to prove that (1) $\Rightarrow$ (2). 

Note that $\operatorname{depth}(M \otimes_R N) \neq n$ by (1). Since $M \otimes_R N$ satisfies $(S_{n+1})$-condition on the punctured spectrum and $n \geq 0$, then $M \otimes_R N$ satisfies at least the $(S_1)$-condition on the punctured spectrum. So, the $R$-modules $M$ and $N$ satisfy the depth formula by Theorem \ref{teo:torindependent} and Corollary \ref{wef}.  Therefore 
\begin{align*}
    \depth (M \otimes_R N) & = \depth M + \depth N - \depth R \\
    & \geq \depth R + n - \depth R =n.
\end{align*}
Thus $\depth (M \otimes_R N)_{\mathfrak{m}} \geq n+1 = \min \lbrace n+1, \dim R_{\mathfrak{m}} \rbrace$ and $M \otimes_R N$ satisfies $(S_{n+1})$-condition.
\end{proof}

\section{Applications to modules of finite $C$-quasi-projective dimension}
In this section, our aim is to establish a more general version of Auslander's depth formula for $C$-quasi-projective dimension, complementing the case proved for $\operatorname{Tor}$-independent modules under certain additional conditions in \cite[Theorem 6.10]{ferraro}. As an application, we provide a new proof of the dependency formula (see \cite[Definition 3.5]{ferraro}) for $C$-quasi-projective dimension (\cite[Corollary 6.15]{ferraro}) based on a direct argument that avoids the reduction to the case $C=R$. More than it, our result (Theorem \ref{depenency}) improves the result by Dey-Ferraro-Gheibi by showing that $q^R(M,N)$ can be computed from the (finitely many) prime ideals in $\operatorname{Ass}(\operatorname{Tor}_{q^R(M,N)}^R(M,N))$.

Before, we will state the recent definition of quasi-projective dimension with respect to a semidualizing module introduced in \cite{ferraro}. 
\begin{newclaim}[\textbf{Semidualizing modules}] A finitely generated $R$-module $C$ is called a \textit{semidualizing} $R$-module if the natural homothety map $R \to \operatorname{Hom}_R(C,C)$ is an isomorphism and $\operatorname{Ext}_R^i(C,C)=0$ for all $i \geq 1$. 
\end{newclaim}
\begin{newclaim}[\textbf{Auslander and Bass classes}]
Let $C$ be a semidualizing $R$-module.
The \textit{Auslander class} $\mathcal{A}_C (R)$ is the class of $R$-modules $M$ satisfying in the following conditions: 
\begin{enumerate}
    \item The natural map $M \rightarrow \operatorname{Hom}_R(C, C \otimes_R M)$ is an isomorphism. 
    \item One has $\operatorname{Tor}_{>0}^R(C,M)=0=\operatorname{Ext}_R^{>0}(C,C \otimes_R M)$.
\end{enumerate}
The \textit{Bass class} $\mathcal{B}_C(R)$ is the class of $R$-modules $M$ satisfying in the following conditions.
\begin{enumerate}
    \item The evaluation map $C \otimes_R \operatorname{Hom}_R(C,M) \rightarrow M$ is an isomorphism.
    \item One has $\operatorname{Ext}_{R}^{>0} (C,M)=0=\operatorname{Tor}_{>0}^R(C,\operatorname{Hom}_R(C,M))$.
\end{enumerate}
\end{newclaim}
The authors refer to the classical reference \cite{sather} for more on semidualizing modules and Auslander-Bass classes.

\begin{definition}(\cite[Definition 4.1]{ferraro})
Let $C$ be a semidualizing $R$-module. An $R$-module $M$ is said to have \textit{finite $C$-quasi-projective dimension} if there exists a bounded complex $P$ of projective $R$-modules such that $P \otimes_R C$ is not acyclic and all the homologies are a finite direct sum of copies of $M$ (or zero). Such a complex $P$ is said to be a $C$-quasi-projective resolution of $M$. The $C$-quasi-projective dimension of $M$ is defined as 
\begin{align*}
\operatorname{\textit{C}-qpd}_R M= \inf \lbrace \sup (P \otimes_R C) - \operatorname{hsup} (P \otimes_R C) \mid P \text{ is a $C$-quasi-projective resolution of $M$} \rbrace.
\end{align*}
The $C$-quasi-projective dimension of the zero module is set to be $- \infty$.
\end{definition}

We remark that when $C=R$, the definition above recovers Definition \ref{defqpd}.  As an application of Theorem \ref{auslanderqpd}, we have the following theorem, which provides an extension of Auslander's original result (Theorem \ref{austordep}) in terms of $C$-quasi-projective dimension. We should mention that the case $q^R(M,N)=0$ in the theorem has been proved in \cite[Theorem 6.10]{ferraro}.

\begin{theorem}\label{depthformulaC}
Let $R$ be a local ring, $M$ and $N$ be non-zero $R$-modules such that $q^R(M,N)< \infty$ and let $C$ be a semidualizing $R$-module. If the following conditions are satisfied 
\begin{enumerate}
    \item $\operatorname{\textit{C}-qpd}_R M < \infty$ and $M \in \mathcal{B}_C(R)$, 
    \item $N \in \mathcal{A}_C(R)$, 
    \item $q^R(M,N)=0$ or $\operatorname{depth}(\operatorname{Tor}_{q^R(M,N)}^{R}(M,N)) \leq 1$,
\end{enumerate}
then $$\operatorname{depth} N + \operatorname{depth} M = \operatorname{depth} R+ \operatorname{depth}(\operatorname{Tor}_{q^R(M,N)}^R(M,N)) -q^R(M,N).$$
\begin{proof} Since $\operatorname{\textit{C}-qpd}_R M< \infty$ and $ M  \in \mathcal{B}_C (R)$, it follows that $\operatorname{qpd}_R (\operatorname{Hom}_R(C,M)) < \infty$, by \cite[Theorem 4.8]{ferraro}. Moreover, since $M \in \mathcal{B}_C(R)$ and $N \in \mathcal{A}_C(R)$, we have $\operatorname{Tor}_i^R(M,N) \cong \operatorname{Tor}_i^R(\operatorname{Hom}_R(C,M), C \otimes_R N)$ for all $i \geq 0$, by \cite[Lemma 3.1.13(c)]{sather}.  Thus, 
\begin{align*}
  q:= q^R(\operatorname{Hom}_R(C,M),C \otimes_R N) = q^R(M,N) < \infty 
\end{align*}
and
\begin{align}\label{isodualizing}
    \operatorname{Tor}_q^R(M,N) \cong \operatorname{Tor}_q^R(\operatorname{Hom}_R(C,M), C \otimes_R N). 
\end{align}
We are now in conditions to apply Corollary \ref{wef} to the pair $\operatorname{Hom}_R(C,M)$ and $C \otimes_R N$. Thus
\begin{align}\label{eq10}
\operatorname{depth} C \otimes_R N + \operatorname{depth} (\operatorname{Hom}_R(C,M))  = \operatorname{depth} R  + \operatorname{depth}\operatorname{Tor}_q^R(\operatorname{Hom}_R(C,M), C \otimes_R N) - q.
\end{align}
By (\ref{isodualizing}), we have $\operatorname{depth}(\operatorname{Tor}_q^R(\operatorname{Hom}_R(C,M), C \otimes_R N))= \operatorname{depth}(\operatorname{Tor}_q^R(M,N))$. Moreover, we have equalities $\operatorname{depth} (\operatorname{Hom}_R(C,M))= \operatorname{depth} M$ and $ \operatorname{depth} C \otimes_R N = \operatorname{depth} N $  by \cite[Lemma 3.9]{aihara} and \cite[Lemma 2.11(i)]{dibaei}, respectively. Using all these equalities and \eqref{eq10}, the result follows. 

\end{proof}
\end{theorem}

We are now almost ready to give a new proof that is different from the one given in \cite[Corollary 6.15]{ferraro} for the dependency formula for $C$-quasi-projective dimension through a direct argument that avoids the reduction to the case $C=R$ and is deduced from Theorem \ref{depthformulaC}. Before, we need to establish the following lemma, that is a version of Corollary \ref{ineqgs} for $C$-quasi-projective dimension.

\begin{lemma}\label{lemaCqpd}
Let $R$ be a local ring, $M$ and $N$ be non-zero $R$-modules such that $q^R(M,N)< \infty$ and let $C$ be a semidualizing $R$-module. If the following conditions are satisfied 
\begin{enumerate}
    \item $\operatorname{\textit{C}-qpd}_R M < \infty$ and $M \in \mathcal{B}_C(R)$, 
    \item $N \in \mathcal{A}_C(R)$,
    \end{enumerate}
then $q^R(M,N) \geq \operatorname{\textit{C}-qpd}_R M - \operatorname{depth} N$.
\begin{proof} 
Since $\operatorname{\textit{C}-qpd}_R M< \infty$ and $ M  \in \mathcal{B}_C (R)$, we must have that $\operatorname{\textit{C}-qpd}_R M= \operatorname{qpd}_R (\operatorname{Hom}_R(C,M)) < \infty$, by \cite[Theorem 4.8]{ferraro}. Moreover, since $M \in \mathcal{B}_C(R)$ and $N \in \mathcal{A}_C(R)$, we have $\operatorname{Tor}_i^R(M,N) \cong \operatorname{Tor}_i^R(\operatorname{Hom}_R(C,M), C \otimes_R N)$ for all $i \geq 0$, by \cite[Lemma 3.1.13(c)]{sather}. That is,
\begin{align*}
  q:= q^R(\operatorname{Hom}_R(C,M),C \otimes_R N) = q^R(M,N) < \infty.
\end{align*}
Thus, we can apply Corollary \ref{ineqgs} to the pair $(\operatorname{Hom}_R(C,M),C \otimes_R N)$, and have the inequality: $q \geq \operatorname{qpd}_R (\operatorname{Hom}_R(C,M)) - \operatorname{depth} (C \otimes_R N)$. Since $\operatorname{\textit{C}-qpd}_R M= \operatorname{qpd}_R( \operatorname{Hom}_R(C,M))$ and $\operatorname{depth} C \otimes_ R N = \operatorname{depth} N$ (see \cite[Lemma 2.11(1)]{dibaei}), we  have the desired inequality. 
\end{proof}
\end{lemma}
The following result gives a new proof of the dependency formula for $C$-quasi-projective dimension \cite[Corollary 6.15]{ferraro} and more than this, it improves the result by pointing out that $q^R(M,N)$ can be computed from the (finitely many) prime ideals in $\operatorname{Ass}_R(\operatorname{Tor}_{q^R(M,N)}^R(M,N))$. 

\begin{theorem}\label{depenency}
Let $R$ be a local ring, $M$ and $N$ be non-zero $R$-modules such that $q^R(M,N)< \infty$ and let $C$ be a semidualizing $R$-module. If the following conditions are satisfied 
\begin{enumerate}
    \item $\operatorname{\textit{C}-qpd}_R M < \infty$ and $M \in \mathcal{B}_C(R)$, 
    \item $N \in \mathcal{A}_C(R)$,
 \end{enumerate}
 then
 \begin{align*}
    q^R(M,N) & = \sup \lbrace \operatorname{depth} R_{\mathfrak{p}} - \operatorname{depth} M_{\mathfrak{p}} - \operatorname{depth} N_{\mathfrak{p}} \mid \mathfrak{p} \in \Spec R \rbrace \textnormal{ (\cite[Corollary 6.15]{ferraro})} \\
    & =  \sup \lbrace \operatorname{depth} R_{\mathfrak{p}} - \operatorname{depth} M_{\mathfrak{p}} - \operatorname{depth} N_{\mathfrak{p}} \mid \mathfrak{p} \in \operatorname{Ass} (\operatorname{Tor}_{q^R(M,N)}^R(M,N)) \rbrace.
\end{align*}
\begin{proof}
Let $\mathfrak{p} \in \Spec R$. If $\mathfrak{p} \notin \operatorname{Supp} M \cap \operatorname{Supp} N$, the inequality we want to prove $q^R(M,N) \geq \operatorname{depth}R_{\mathfrak{p}}-\depth_{R_{\mathfrak{p}}} M_{\mathfrak{p}}-\depth_{R_{\mathfrak{p}}} N_{\mathfrak{p}}$ is trivial. So, in order to prove this inequality, we may assume that $\mathfrak{p} \in \operatorname{Supp} M \cap \operatorname{Supp} N$. Note that $q^{R_{\mathfrak{p}}}(M_{\mathfrak{p}},N_{\mathfrak{p}}) \leq q^R(M,N) < \infty$ and that $C_{\mathfrak{p}}$ is a semidualizing module over $R_{\mathfrak{p}}$. Moreover, one can see that $\operatorname{\textit{C}_{\mathfrak{p}}-qpd}_{R_{\mathfrak{p}}} M_{\mathfrak{p}} \leq \operatorname{\textit{C}-qpd}_R M < \infty$, by \cite[Theorem 5.7]{ferraro}. Since $M_{\mathfrak{p}} \in \mathcal{B}_{C_{\mathfrak{p}}}(R_{\mathfrak{p}})$ and $N_{\mathfrak{p}} \in \mathcal{A}_{C_{\mathfrak{p}}}(R_{\mathfrak{p}})$, Lemma \ref{lemaCqpd} applies to $(M_{\mathfrak{p}},N_{\mathfrak{p}})$. This yields:
\begin{align}\label{inequalities}
    q^R(M,N) \geq q^{R_{\mathfrak{p}}}(M_{\mathfrak{p}},N_{\mathfrak{p}})
    \geq \operatorname{\textit{C}_{\mathfrak{p}}-qpd}_{R_{\mathfrak{p}}} M_{\mathfrak{p}} - \operatorname{depth}_{R_{\mathfrak{p}}} N_{\mathfrak{p}}.
\end{align}
So, using the Auslander-Buchsbaum formula for $C$-quasi-projective dimension (\cite[Theorem 5.10]{ferraro}), (\ref{inequalities}) implies that $q^R(M,N) \geq \operatorname{depth} R_{\mathfrak{p}} - \operatorname{depth}_{R_{\mathfrak{p}}} M_{\mathfrak{p}} - \operatorname{depth}_{R_{\mathfrak{p}}} N_{\mathfrak{p}}$ for all $\mathfrak{p} \in \operatorname{Spec}R $. Then 
\begin{align*}
q^R(M,N) & \geq \sup \lbrace \operatorname{depth} R_{\mathfrak{p}} - \operatorname{depth} M_{\mathfrak{p}} - \operatorname{depth} N_{\mathfrak{p}} \mid \mathfrak{p} \in \Spec R \rbrace \\
& \geq   \sup \lbrace \operatorname{depth} R_{\mathfrak{p}} - \operatorname{depth} M_{\mathfrak{p}} - \operatorname{depth} N_{\mathfrak{p}} \mid \mathfrak{p} \in \operatorname{Ass} (\operatorname{Tor}_{q^R(M,N)}^R(M,N)) \rbrace.
\end{align*}

Now, let $\mathfrak{p} \in  \operatorname{Ass} (\operatorname{Tor}_{q^R(M,N)}^R(M,N))$. Then $\operatorname{depth}_{R_{\mathfrak{p}}} (\operatorname{Tor}_{q^R(M,N)}^{R_{\mathfrak{p}}}(M_{\mathfrak{p}},N_{\mathfrak{p}})) =0$, and, because $\operatorname{depth} 0=\infty$, it follows that $\operatorname{Tor}_{q^R(M,N)}^{R_{\mathfrak{p}}}(M_{\mathfrak{p}},N_{\mathfrak{p}}) \neq 0$ and thus $q^{R_{\mathfrak{p}}}(M_{\mathfrak{p}},N_{\mathfrak{p}}) =q^R(M,N)$. Given that $\operatorname{depth} (\operatorname{Tor}_{q^R(M,N)}^{R_{\mathfrak{p}}}(M_{\mathfrak{p}},N_{\mathfrak{p}})) =0$, applying Theorem \ref{depthformulaC} to the pair $M_{\mathfrak{p}}$ and $N_{\mathfrak{p}}$ we obtain the equality:
$$q^R(M,N) = \operatorname{depth} R_{\mathfrak{p}} - \operatorname{depth} M_{\mathfrak{p}} - \operatorname{depth} N_{\mathfrak{p}}$$ for all  $\mathfrak{p} \in \operatorname{Ass} (\operatorname{Tor}_{q^R(M,N)}^R(M,N))$.
As this equality is true for each associated prime ideal of $\operatorname{Tor}_{q^R(M,N)}^R(M,N)$, we deduce that 
$$q^R(M,N)=\sup \lbrace \operatorname{depth} R_{\mathfrak{p}} - \operatorname{depth} M_{\mathfrak{p}} - \operatorname{depth} N_{\mathfrak{p}} \mid \mathfrak{p} \in \operatorname{Ass} (\operatorname{Tor}_{q^R(M,N)}^R(M,N)) \rbrace.$$
and the result follows. 
\end{proof}
\end{theorem}
\begin{remark}
When $C=R $, one obtains from Theorem \ref{depenency} an improved version of \cite[Theorem 3.6]{ferraro} showing that $q^R (M,N)$ can be computed from the (finitely many) prime ideals in $\operatorname{Ass}_R(\operatorname{Tor}_{q^R(M,N)}^R(M,N))$ when $M$ is of finite quasi-projective dimension. 
\end{remark}

\begin{agra}
Part of this work was developed when the third author was visiting the University of Nebraska-Lincoln. He is grateful for the hospitality. This work was completed while the second author was visiting the Department of Mathematics at University of Texas at Arlington, and he gratefully acknowledges their hospitality. The authors are grateful to the anonymous referee for a
careful reading and for valuable comments.
\end{agra} 

\noindent
\textbf{Disclosure statement.} No potential conflict of interest was reported by the authors.

\begin{funding}
    The first author was supported by grant 2019/21181-0, S\~ao Paulo Research Foundation (FAPESP) and  grant 304851/2025-6, CNPq-Brazil. The second author was supported by grants 2022/12114-0 and 2024/17809-1, S\~ao Paulo Research Foundation (FAPESP). The third author was supported by grants 2022/03372-5 and 2023/15733-5, S\~ao Paulo Research Foundation (FAPESP).
\end{funding}

\end{document}